\documentclass[11pt,a4paper]{article}
\usepackage[fleqn]{amsmath}
\usepackage{amsthm}
\usepackage{amssymb}

\newtheorem{proposition}{Proposition}
\newtheorem{corollary}[proposition]{Corollary}
\newtheorem{theorem}[proposition]{Theorem}
\newtheorem{lemma}[proposition]{Lemma}

\theoremstyle{definition}
\newtheorem{definition}{Definition}
\newtheorem{example}{Example}

\newcommand{\goth}[1]{\mathfrak{#1}}

\newcommand{\schl}{[\![}
\newcommand{\schr}{]\!]}

\DeclareMathOperator{\coad}{coad} \DeclareMathOperator{\Ad}{Ad}
\DeclareMathOperator{\id}{id} \DeclareMathOperator{\Alt}{Alt}
\DeclareMathOperator{\Lie}{Lie} 
\DeclareMathOperator{\sign}{sign} \DeclareMathOperator{\CYB}{CYB}
\DeclareMathOperator{\ad}{ad}

\begin{document}
\title{On quasi-Poisson homogeneous spaces of quasi-Poisson Lie groups}
\author{Eugene Karolinsky, Kolya Muzykin}
\date{}
\maketitle
\section{Introduction}
The notion of Poisson Lie group and its infinitesimal counterpart,
Lie bialgebra, was introduced by Drinfeld \cite{D_PL}. Later it
was explained that these objects are quasiclassical limits of Hopf
QUE algebras. In \cite{D_quasihopf_alg} the more general objects,
quasi-Hopf QUE algebras, were introduced along with their
quasiclassical limits, Lie quasi-bialgebras. The corresponding
geometric objects, quasi-Poisson Lie groups, were first studied by
Kosmann-Schwarzbach \cite{YKS_qbialg_qpoiss_fr}.

It is well known that Lie bialgebra structures on $\goth{g}$ are
in a natural 1-1 cor\-res\-pon\-dence with Lie algebra structures
on $\mathcal{D}(\goth{g})=\goth{g}\oplus\goth{g}^*$ such that
$\goth{g}$ and $\goth{g}^*$ are subalgebras in
$\mathcal{D}(\goth{g})$ and the natural bilinear form on
$\mathcal{D}(\goth{g})$ is invariant. Respectively, in order to
get a Lie quasi-bialgebra structure on $\goth{g}$, one should drop
the condition that $\goth{g}^*$ is a subalgebra in
$\mathcal{D}(\goth{g})$.

Along with (quasi-)Poisson Lie groups it is natural to study their
\makebox{(quasi-)}\linebreak[1]Poisson actions
\cite{AAYKS_mpairs_mmaps, AAYKSEM} and, in particular,
(quasi-)Poisson homogeneous spaces. Drinfeld in
\cite{D_poiss_hom_spaces} presented an approach to the
classification of Poisson homogeneous spaces. Namely, he showed
that if $G$ is a Poisson Lie group, $\goth{g}$ is the
corresponding Lie bialgebra, then the isomorphism classes of
Poisson homogeneous $G$-spaces are essentially in a 1-1
correspondence with the $G$-orbits of Lagrangian subalgebras in
$\mathcal{D}(\goth{g})$.

The main goal of this paper is to generalize this result to the
quasi-Poisson case (see Theorem \ref{mainth}). We also study the
behavior of quasi-Poisson homogeneous spaces under twisting. Some
examples showing the technique of Lagrangian subalgebras are also
provided.

It also turns out that quasi-Poisson homogeneous spaces, as well
as Poisson ones, are related to solutions of the classical
dynamical Yang-Baxter equation (see \cite{KS_dynam_rmatr,
Lu_CDYBE} for the Poisson case). This topic will be discussed in a
forthcoming paper.
\subsection*{Acknowledgements} The authors are grateful to
Alexander Stolin for useful discussions on the topic of the paper.
\section{Preliminaries}
\subsection{Notation}
We will use the following normalization of the wedge product of
multivector fields on a smooth manifold. If $v$ is a $m$-vector
field, $w$ is a $n$-vector field, then
\begin{equation*}
v\wedge w=\frac{1}{n!m!}\Alt(v\otimes w),
\end{equation*}
where
\begin{equation*}
\Alt(x_1\otimes x_2\otimes\dotsm\otimes x_k)=\sum_{\sigma\in
S_k}\sign(\sigma)x_{\sigma(1)}\otimes x_{\sigma(2)}\otimes
\dotsm\otimes x_{\sigma(k)}.
\end{equation*}

We will denote by $[\![\;,\;]\!]$ the Schouten bracket of
multivector fields (see, e.g., \cite{AAYKSEM}).

Let $G$ be a Lie group, $\goth{g}=\Lie G$ its Lie algebra. For any
$v\in\bigwedge^\bullet\goth{g}$ denote by $v^\lambda$ (resp.
$v^\rho$) the left (resp. right) invariant multivector field that
corresponds to $v$, i.e., $v^\lambda(g)=(l_g)_*v$,
$v^\rho(g)=(r_g)_*v$ for all $g\in G$, where $l_g$ (resp. $r_g$)
is the left (resp. right) translation by $g$.

Suppose that $G$ acts smoothly on a smooth manifold $X$. Then for
any $v\in\goth{g}$ we denote by $v_X$ the corresponding vector
field on $X$, i.e.,
\begin{equation*}
(v_Xf)(x)=\frac{d}{dt}\Bigl|_{t=0}f(\exp tv\cdot x )
\end{equation*}
for any $x\in X$. Similarly, for $v\in\bigwedge^\bullet\goth{g}$
one can define the multivector field $v_X$. For any $x\in X$
consider the map $\rho_x:G\to X$, $\rho_x(g)=g\cdot x$. Then
$(\rho_x)_*v=v_X(x)$ for $v\in\goth{g}$.

For any point $x\in X$ we denote by $H_x=\{g\in G\mid g\cdot
x=x\}$ its stabilizer. Let $\goth{h}_x=\Lie H_x\subset\goth{g}$.

Suppose now that $X$ is a homogeneous $G$-space. In this case we
will identify $T_xX$ with $\goth{g/h}_x$ for all $x\in X$. Fix
$x\in X$ and for any $f\in C^\infty X$ define $f^G\in C^\infty G$
by the formula $f^G(g)=(f\circ\rho_x)(g)=f(g\cdot x)$. Note that
the mapping $f\mapsto f^G$ is an isomorphism between the spaces of
smooth functions on $X$ and right $H_x$-invariant smooth functions
on $G$.

\subsection{Quasi-Poisson Lie groups and quasi-Poisson actions}
Following \cite{AAYKS_mpairs_mmaps}, we define the notion of
quasi-Poisson Lie group and the notion of quasi-Poisson action.
\begin{definition}
Let $G$ be a Lie group, $\goth{g}$ its Lie algebra, $P_G$ a
bivector field on $G$, and $\varphi\in\bigwedge^3\goth{g}$. A
triple $(G,P_G,\varphi)$ is called a \emph{quasi-Poisson Lie
group} if
\begin{eqnarray}\label{multcond}
&&\hspace{-1cm}\text{$P_G$ is \emph{multiplicative}, i.e., }
P_G(gg')=(l_g)_*P_G(g')+(r_{g'})_*P_G(g), \\
\label{quasipoiss} &&\hspace{-1cm}\frac{1}{2}\schl
P_G,P_G\schr=\varphi^\rho-\varphi^\lambda,\\
&&\hspace{-1cm}\schl P_G,\varphi^\rho\schr=0.
\end{eqnarray}
\end{definition}

The notion of Poisson Lie group is a special case of the notion of
quasi-Poisson Lie group. Namely, for any Poisson Lie group
$(G,P_G)$ the triple $(G,P_G,0)$ is a quasi-Poisson Lie group.

Consider the mapping $\eta:G\rightarrow\goth{g}\wedge\goth{g}$
defined by $\eta(g)=(r_g^{-1})_*P_G(g)$. It is a
$\goth{g}\wedge\goth{g}$-valued 1-cocycle of $G$ with respect to
the adjoint action of $G$ on $\goth{g}\wedge\goth{g}$, i.e.,
\begin{equation*}
\eta(g_1g_2)=\eta(g_1)+\Ad_{g_1}\eta(g_2).
\end{equation*}
Here $\Ad_g(x\otimes y)=(\Ad_gx)\otimes(\Ad_gy)$. The cocyclicity
of $\eta$ is equivalent to the multiplicativity condition
\eqref{multcond}.

Consider
$\delta=d_e\eta:\goth{g}\rightarrow\goth{g}\wedge\goth{g}$. It is
a 1-cocycle of $\goth{g}$ with respect to the adjoint action of
$\goth{g}$ on $\goth{g}\wedge\goth{g}$, i.e.,
\begin{equation*}
\delta([x,y])=\ad_x\delta(y)-\ad_y\delta(x),
\end{equation*}
where $\ad_x(y\otimes z)= [x\otimes 1+1\otimes x,y\otimes
z]=\ad_xy\otimes z+y\otimes\ad_x z$.
\begin{definition}
Suppose $(G,P_G,\varphi)$ is a quasi-Poisson Lie group, $G$ acts
smooth\-ly on a smooth manifold $X$, $P_X$ is a bivector field on
$X$. The action of $G$ on $X$ is called \emph{quasi-Poisson} if
\begin{eqnarray}\label{quasipoissonmul}
&&\hspace{-5.35cm}P_X(gx)=(l_g)_*P_X(x)+(\rho_x)_*P_G(g), \\
\label{quasipoissonjac} &&\hspace{-5.35cm}\frac{1}{2}\schl
P_X,P_X\schr=\varphi_X
\end{eqnarray}
(here $l_g$ denotes the mapping $x\mapsto g\cdot x$).
\end{definition}

Let us consider the case $\varphi=0$, i.e., $G$ is a Poisson Lie
group. Then the condition \eqref{quasipoissonjac} means that $X$
is a Poisson manifold, and from \eqref{quasipoissonmul} it follows
that the action of $G$ on $X$ is Poisson.

\begin{definition}
Suppose that $(G,P_G,\varphi)$ is a quasi-Poisson group, $G$ acts
smoothly on a manifold $X$ equipped with a bivector field $P_X$,
and this action is quasi-Poisson. We call $X$ a
\emph{quasi-Poisson homogeneous $G$-space} if the action of $G$ on
$X$ is transitive.
\end{definition}
\begin{lemma}
Suppose that $(G,P_G,\varphi)$ is a quasi-Poisson group, $X$ is a
ho\-mo\-ge\-neous $G$-space, $P_X$ is a bivector field on $X$.
Then the condition \eqref{quasipoissonmul} is equivalent to
\begin{equation}\label{quasipoissonmulhom}
P_X(gx)=\Ad_gP_X(x)+\overline{\eta(g)},
\end{equation}
where
$\Ad_g:\bigwedge^2(\goth{g/h}_x)\to\bigwedge^2(\goth{g/h}_{gx})$
is the isomorphism of the vector spaces induced by the
automorphism $\Ad_g:\goth{g}\to\goth{g}$, and $\overline{\eta(g)}$
is the image of $\eta(g)$ in $\bigwedge^2(\goth{g/h}_{gx})$. \qed
\end{lemma}

\subsection{Lie quasi-bialgebras}
Recall that a Poisson Lie structure on a Lie group $G$ induces the
structure of a Lie bialgebra on the Lie algebra $\goth{g}=\Lie G$.
A quasi-Poisson structure on a Lie group $G$ induces a similar
structure on $\goth{g}$. We follow \cite{D_quasihopf_alg} in
defining the notion of Lie quasi-bialgebra.
\begin{definition}
Let $\goth{g}$ be a Lie algebra, $\delta$ a
$\goth{g}\wedge\goth{g}$-valued 1-cocycle of $\goth{g}$, and
$\varphi\in\bigwedge^3\goth{g}$. A triple
$(\goth{g},\delta,\varphi)$ is called a \emph{Lie quasi-bialgebra}
if
\begin{eqnarray}\label{quasicojac}
&&\hspace{-4.6cm}\frac{1}{2}\Alt(\delta\otimes\id)\delta(x)=\ad_x\varphi\quad\text{for
any $x\in\goth{g}$}, \\
\label{infpentagon}
&&\hspace{-4.6cm}\Alt(\delta\otimes\id\otimes\id)\varphi=0,
\end{eqnarray}
where $\ad_x(a\otimes b\otimes c)=[x\otimes 1\otimes 1+1\otimes
x\otimes 1+1\otimes 1\otimes x,a\otimes b\otimes c]$.
\end{definition}

The equation \eqref{quasicojac} is called the quasi co-Jacobi
identity.

If we set $\varphi=0$, then the notion of Lie quasi-bialgebra
coincides with the notion of Lie bialgebra. In this case the
equation \eqref{quasicojac} becomes the ordinary co-Jacoby
identity, and the condition \eqref{infpentagon} is obviously
satisfied.

For any quasi-Poisson Lie group $(G,P_G,\varphi)$ there exists a
Lie quasi-bi\-al\-geb\-ra structure on $\goth{g}$ given by the
1-cocycle $\delta=d_e\eta$ and $\varphi$. Conversely, to any Lie
quasi-bialgebra there corresponds a unique connected and simply
connected quasi-Poisson Lie group (see \cite{YKS_qbialg_qpoiss}).

Given any linear map
$\delta:\goth{g}\to\goth{g}\wedge\goth{g}\subset\goth{g}\otimes\goth{g}$
we can define the skew-symmetric bilinear operation on
$\goth{g}^*$: for all $l,m\in\goth{g}^*$ set
$[l,m]_{\delta}=\delta^*(l\otimes m)$.

Recall that for any Lie quasi-bialgebra
$(\goth{g},\delta,\varphi)$ one can construct the so-called
\emph{double Lie algebra} $\mathcal{D}(\goth{g})$ (see
\cite{MBYKS_djacbialgebra}):

let $\mathcal{D}(\goth{g})=\goth{g}\oplus\goth{g}^*$ as a vector
space;

define the bilinear operation $[\;,\;]_{\mathcal{D}(\goth{g})}$ on
$\mathcal{D}(\goth{g})$ by the following conditions:
\begin{enumerate}
\item $[a,b]_{\mathcal{D}(\goth{g})}=[a,b]\quad$  for $a,b\in\goth{g}$;
\item $[l,m]_{\mathcal{D}(\goth{g})}=[l,m]_{\delta}-(l\otimes m\otimes\id)\varphi\quad$ for
$l,m\in\goth{g}^*$;
\item $[a,l]_{\mathcal{D}(\goth{g})}=\coad_al-\coad_la\quad$ for $a\in\goth{g},
l\in\goth{g}^*$.
\end{enumerate}
where $\coad_l:\goth{g}\to\goth{g}$ is defined by
$\bigl<\coad_la,m\bigr>=-\bigl<[l,m]_{\delta},a\bigr>=-\bigl<l\otimes
m,\delta(a)\bigr>$, and $\coad_a:\goth{g}^*\to\goth{g}^*$ is
defined by $\bigl<\coad_al,b\bigr>=-\bigl<l,[a,b]\bigr>$. Here and
below $\bigl<\;,\;\bigr>$ denotes the standard pairing between
$\goth{g}$ and $\goth{g}^*$.

We denote by $Q(\;,\;)$ the following invariant symmetric bilinear
form on $\mathcal{D}(\goth{g})$:
\begin{equation*}
Q(a+l,b+m)=\bigl<l,b\bigr>+\bigl<m,a\bigr>.
\end{equation*}

Suppose $G$ is a quasi-Poisson Lie group, $\goth{g}$ is the
corresponding Lie quasi-bialgebra, $\mathcal{D}(\goth{g})$ is its
double Lie algebra. Then the adjoint action of $G$ on $\goth{g}$
can be extended to the action of $G$ on $\mathcal{D}(\goth{g})$
defined by
\begin{equation*}
g\cdot(a+l)=\Ad_ga+(l'\otimes\id)\eta(g)+l',
\end{equation*}
where $l'=(\Ad_g^{-1})^*l$. The differential of this action is the
adjoint action of $\goth{g}$ on $\mathcal{D}(\goth{g})$.
\section{Main results}
In \cite{D_poiss_hom_spaces} the characterization of all Poisson
homogeneous structures on a given ho\-mo\-ge\-neous $G$-space in
terms of Lagrangian subalgebras in $\mathcal{D}(\goth{g})$ is
presented. We generalize this result to the quasi-Poisson case.

Suppose $G$ is a quasi-Poisson Lie group, $X$ is a quasi-Poisson
ho\-mo\-ge\-neous $G$-space. Recall that we identify $T_xX$ and
$\goth{g/h}_x$ for all $x\in X$. For any $x\in X$ define
\begin{equation*}
L_x=\left\{a+l\mid a\in\goth{g},\
l\in(\goth{g/h}_x)^*=\goth{h}_x^\perp\subset\goth{g}^*,\
(l\otimes\id)P_X(x)=\overline{a}\right\},
\end{equation*}
where $\overline{a}$ is the image of $a$ in $\goth{g/h}_x$.

\begin{lemma}
$L_x$ is Lagrangian (that is, maximal isotropic) subspace in
$\mathcal{D}(\goth{g})$, and $L_x\cap\goth{g}=\goth{h}_x$. \qed
\end{lemma}
Denote by $\Lambda$ the set of all Lagrangian subalgebras in
$\mathcal{D}(\goth{g})$.
\begin{theorem}\label{mainth}
Suppose $(G,P_G,\varphi)$ is a quasi-Poisson Lie
group, $(X,P_X)$ is a quasi-Poisson homogeneous $G$-space. Then
the following statements hold:
\begin{enumerate}
\item $L_x$ is a subalgebra in $\mathcal{D}(\goth{g})$
for all $x\in X$;
\item $L_{gx}=g\cdot L_x$;
\item Thus we get a bijection between the set of all $G$-quasi-Poisson
structures on $X$ and the set of $G$-equivariant maps $x\mapsto
L_x$ from $X$ to $\Lambda$ such that $L_x\cap\goth{g}=\goth{h}_x$
for all $x\in X$.
\end{enumerate}
\end{theorem}
\begin{corollary}
There is a bijection between the set of all isomorphism classes of
quasi-Poisson homogeneous $G$-spaces and the set of $G$-conjugacy
classes of pairs $(L,H)$, where $L\subset\mathcal{D}(\goth{g})$ is
a Lagrangian subalgebra, $H$ is a closed subgroup in $G_L=\{g\in
G\mid g\cdot L=L\}$, and $L\cap\goth{g}=\Lie H$. \qed
\end{corollary}
The rest of this section is devoted to the proof of Theorem
\ref{mainth}. We start with a technical lemma.

\newcommand{\locderiv}[1]{\partial_{#1}}

\begin{lemma}\label{schoutensquarelem}
Let $P$ be a bivector field on a smooth manifold $X$. Define
$\{f_1,f_2\}=P(df_1,df_2)$ for all $f_1,f_2\in C^\infty X$. Then
\begin{equation}
\oint\{\{f_1,f_2\},f_3\}=-\frac{1}{2}\schl
P,P\schr(df_1,df_2,df_3),
\end{equation}
where $\oint$ denotes the sum over all cyclic permutations of
$f_1, f_2, f_3$
\end{lemma}
\begin{proof}
Straightforward computation.
\end{proof}

\begin{lemma}\label{equivariancelem}
$L_{gx}=g\cdot L_x$ iff \eqref{quasipoissonmulhom} holds.
\end{lemma}
\begin{proof}
By definition,
\begin{equation*}
L_x=\left\{a+l\mid a\in\goth{g}, l\in(\goth{g/h}_x)^*,
(l\otimes\id)P_X(x)=\overline{a}\right\},
\end{equation*}
\begin{equation*}
L_{gx}=\left\{a'+l'\mid a'\in\goth{g}, l'\in(\goth{g/h}_{gx})^*,
(l'\otimes\id)P_X(gx)=\overline{a'}\right\}.
\end{equation*}

It is enough to check that
\begin{equation*}
g\cdot L_x=\left\{a'+l'\mid a'\in\goth{g},
l'\in(\goth{g/h}_{gx})^*,
(l'\otimes\id)\bigl(\Ad_gP_X(x)+\overline{\eta(g)}\bigr)=\overline{a'}\right\}.
\end{equation*}

Consider $a'+l'=g\cdot(a+l), a\in\goth{g}, l\in(\goth{g/h}_x)^*$,
that is,
\begin{equation*}
l'=(\Ad_g^{-1})^*l,\; a'=\Ad_ga+(l'\otimes\id)\eta(g).
\end{equation*}

We have
\begin{gather*}
(l'\otimes\id)\bigl(\Ad_gP_X(x)+\overline{\eta(g)}\bigr)=\\
(l\otimes\id)(\Ad_g^{-1}\otimes\id)(\Ad_g\otimes\Ad_g)P_X(x)+
(l'\otimes\id)\overline{\eta(g)}=\\
\Ad_g(l\otimes\id)P_X(x)+(l'\otimes\id)\overline{\eta(g)}.
\end{gather*}

So
$(l'\otimes\id)\bigl(\Ad_gP_X(x)+\overline{\eta(g)}\bigr)=\overline{a'}$
if and only if $a+l\in L_x$. This proves the required equality.
\end{proof}

Now we are heading for the first statement of the theorem.

Let $e_i$ form a basis in $\goth{g}$, $\locderiv{i}$ (resp.
$\locderiv{i}'$) be the right (resp. left) invariant vector field
on $G$ that corresponds to $e_i$.

Suppose $\eta(g)=\eta^{ij}(g)e_i\wedge e_j$. Then
$P_G=\eta^{ij}\locderiv{i}\wedge\locderiv{j}$. Choose any
$r\in\bigwedge^2\goth{g}$ such that the image of $r$ in
$\bigwedge^2(\goth{g/h}_x)$ equals $P_X(x)$. Define
\begin{equation*}
\CYB(r)=[r^{12},r^{13}]+[r^{12},r^{23}]+[r^{13},r^{23}].
\end{equation*}

\begin{lemma}\label{phiplusetc_vanishlem}
Assume that \eqref{quasipoissonmul} holds. Then the image of
\begin{equation*}
\varphi-\CYB(r)+\frac{1}{2}\Alt(\delta\otimes\id)(r)
\end{equation*}
in $\bigwedge^3(\goth{g/h}_x)$ vanishes iff
\eqref{quasipoissonjac} holds.
\end{lemma}

\begin{proof}
From \eqref{quasipoissonmul} it follows that
\begin{gather*}
P_X(gx)(d_{gx}f_1,d_{gx}f_2)=
\left((l_g)_*P_X(x)+(\rho_x)_*P_G(g)\right)(d_{gx}f_1,d_{gx}f_2)=\\
P_X(x)(d_x(f_1\circ l_g),d_x(f_2\circ l_g))+
P_G(g)(d_g(f_1\circ\rho_x),d_g(f_2\circ\rho_x))=\\ r(d_e(f_1\circ
l_g)^G,d_e(f_2\circ l_g)^G)+P_G(g)(d_gf_1^G,d_gf_2^G)=\\
(r^\lambda(g)+P_G(g))(d_gf_1^G,d_gf_2^G).
\end{gather*}

For any $f_1,f_2\in C^\infty G$ define the bracket
\begin{equation*}
\{f_1,f_2\}(g)=(r^\lambda(g)+P_G(g))(d_gf_1,d_gf_2).
\end{equation*}

Using Lemma \ref{schoutensquarelem} we see that
\begin{gather*}
\oint\{\{f_1,f_2\}_X,f_3\}_X(g\cdot x)=
\oint\{\{f_1^G,f_2^G\},f_3^G\}(g)=\\
-\frac{1}{2}\schl
P_G+r^\lambda,P_G+r^\lambda\schr(df_1^G,df_2^G,df_3^G)(g).
\end{gather*}

\begin{lemma}
$\schl P_G+r^\lambda,P_G+r^\lambda\schr=
2\left(\varphi^\rho-\varphi^\lambda+\CYB(r)^\lambda
-\frac{1}{2}\Alt(\delta\otimes\id)(r)^\lambda\right)$
\end{lemma}
\begin{proof}
Using the graded anticommutativity of Schouten bracket, we get
\begin{equation*}
\schl P_G+r^\lambda,P_G+r^\lambda\schr= \schl P_G,P_G\schr+2\schl
P_G,r^\lambda\schr+\schl r^\lambda,r^\lambda\schr.
\end{equation*}

From \eqref{quasipoiss} it follows that
\begin{equation*}
\schl P_G,P_G\schr=2(\varphi^\rho-\varphi^\lambda).
\end{equation*}

We will calculate the rest of the terms on the right hand side
using coordinates. Let $r=r^{ij}e_i\wedge e_j$. Then
$r^\lambda=r^{ij}\locderiv{i}'\wedge\locderiv{j}'$, and
\begin{gather*}\schl r^\lambda,r^\lambda\schr=-4r^{\mu\nu}r^{ij}\schl
\locderiv{\mu}',\locderiv{i}'\schr
\wedge\locderiv{j}'\wedge\locderiv{\nu}'=\\
-4r^{\mu\nu}r^{ij}\Alt\left(\schl
\locderiv{\mu}',\locderiv{i}'\schr
\otimes\locderiv{j}'\otimes\locderiv{\nu}'\right)=
-\Alt\left([r^{13},r^{12}]\right)^\lambda=2\CYB(r)^\lambda.
\end{gather*}

Now we prove that $\schl P_G,r^\lambda\schr=
-\frac{1}{2}\Alt(\delta\otimes\id)(r)^\lambda$. We have
\begin{gather*}
\schl P_G,r^\lambda\schr=
\schl \eta^{\mu\nu}\locderiv{\mu}\wedge\locderiv{\nu},r^{ij}\locderiv{i}'\wedge\locderiv{j}'\schr=\\
r^{ij}\left(\schl\locderiv{i}',\eta^{\mu\nu}\locderiv{\mu}\schr\wedge\locderiv{j}'\wedge\locderiv{\nu}-
\schl\locderiv{j}',\eta^{\mu\nu}\locderiv{\mu}\schr\wedge\locderiv{i}'\wedge\locderiv{\nu}\right)=\\
2r^{ij}\schl\locderiv{i}',\eta^{\mu\nu}\locderiv{\mu}\schr\wedge\locderiv{j}'\wedge\locderiv{\nu}=
-2r^{ij}(\locderiv{i}'\eta^{\mu\nu})\locderiv{\mu}\wedge\locderiv{\nu}\wedge\locderiv{j}'.
\end{gather*}
Using the cocyclicity of $\eta$, we get
\begin{gather*}
\locderiv{i}'\eta^{\mu\nu}(g)e_\mu\wedge e_\nu=
\frac{d}{dt}\Bigl|_{t=0}\eta^{\mu\nu}(g\exp te_i)e_\mu\wedge
e_\nu=\\
\frac{d}{dt}\Bigl|_{t=0}\left(\eta^{\mu\nu}(g)e_\mu\wedge e_\nu+
\Ad_g(\eta^{\mu\nu}(\exp te_i)e_\mu\wedge e_\nu)\right)=\\
\frac{d}{dt}\Bigl|_{t=0}\eta^{kl}(\exp te_i)
(\Ad_g)_k^\mu(\Ad_g)_l^\nu e_\mu\wedge e_\nu=\\
\locderiv{i}'\eta^{kl}(e)(\Ad_g)_k^\mu(\Ad_g)_l^\nu e_\mu\wedge
e_\nu,
\end{gather*}
where $\Ad_ge_k=(\Ad_g)_k^\mu e_\mu$. So,
$\locderiv{i}'\eta^{\mu\nu}(g)=
\locderiv{i}'\eta^{kl}(e)(\Ad_g)_k^\mu(\Ad_g)_l^\nu.$

Continuing our calculations, we have
\begin{gather*}
\schl P_G,r^\lambda\schr(g)=
-2r^{ij}(\locderiv{i}'\eta^{\mu\nu})(g)\locderiv{\mu}(g)
\wedge\locderiv{\nu}(g)\wedge\locderiv{j}'(g)=\\
-2r^{ij}\locderiv{i}'\eta^{kl}(e)(\Ad_g)_k^\mu(\Ad_g)_l^\nu\locderiv{\mu}(g)
\wedge\locderiv{\nu}(g)\wedge\locderiv{j}'(g)=\\
%
%
-2r^{ij}\locderiv{i}'\eta^{\mu\nu}(e)\locderiv{\mu}'(g)\wedge
\locderiv{\nu}'(g)\wedge\locderiv{j}'(g)=\\
-2r^{ij}\locderiv{i}'\eta^{\mu\nu}(e)\Alt\left(\locderiv{\mu}'(g)\otimes
\locderiv{\nu}'(g)\otimes\locderiv{j}'(g)\right)=\\
-r^{ij}\Alt\left((d_e\eta(e_i))^\lambda(g)\otimes\locderiv{j}'(g)\right)=
-r^{ij}\Alt\left(\delta(e_i)\otimes e_j\right)^\lambda(g)=\\
-\frac{1}{2}\left(\Alt(\delta\otimes\id)r\right)^\lambda(g).
\end{gather*}
\end{proof}

Now we finish the proof of Lemma \ref{phiplusetc_vanishlem}. From
the definition of a quasi-Poisson action it follows that
\begin{gather*}
\oint\{\{f_1,f_2\}_X,f_3\}_X(g\cdot x)=
-\varphi_X(df_1,df_2,df_3)(g\cdot x)=\\
-\varphi^\rho(df_1^G,df_2^G,df_3^G)(g).
\end{gather*}
It means that for all $f_1,f_2,f_3\in C^\infty X$ we have
\begin{equation*}
\left(\varphi-\CYB(r)+\frac{1}{2}\Alt(\delta\otimes\id)r\right)^\lambda(df_1^G,df_2^G,df_3^G)=0.
\end{equation*}
Consequently, for all $l,m,n\in\goth{h}_x^\perp$ we get
\begin{equation*}
\bigl<\varphi-\CYB(r)+\frac{1}{2}\Alt(\delta\otimes\id)r,l\otimes
m\otimes n\bigr>=0,
\end{equation*}
which proves the statement of the lemma.
\end{proof}
\begin{lemma}\label{subalglemma}
Assume that \eqref{quasipoissonmul} holds. Then $L_x$ is a
subalgebra in $\mathcal{D}(\goth{g})$ if and only if the image of
the tensor $\varphi+\frac{1}{2}\Alt(\delta\otimes\id)(r)-\CYB(r)$
in $\bigwedge^3(\goth{g}/\goth{h}_x)$ vanishes.
\end{lemma}
\begin{proof}
Consider the mapping $R:\goth{g}^*\to\goth{g}$ that corresponds to
$r\in\bigwedge^2\goth{g}$:
\begin{equation*}
R(l)=(l\otimes\id)r=\sum_il(r_i')r_i'',
\end{equation*}
where $r=\sum_ir_i'\otimes r_i''$.

Then
\begin{gather*}
L_x=\left\{a+l\mid a\in\goth{g}, l\in(\goth{g/h}_x)^*,
(l\otimes\id)\overline{r}=\overline{a}\right\}=\\ \left\{a+l\mid
a\in\goth{g},l\in\goth{h}_x^\perp,\overline{R(l)}=\overline{a}\right\}=
\left\{l+R(l)\mid l\in\goth{h}_x^\perp\right\}+\goth{h}_x.
\end{gather*}

From Lemma \ref{equivariancelem} it follows that $h\cdot
L_x=L_{hx}=L_x$ for any $h\in H_x$. Consequently, for all
$a\in\goth{h}_x$ we have $\ad_a(L_x)\subset L_x$. So $L_x$ is a
Lie subalgebra in $\mathcal{D}(\goth{g})$ if and only if
$[l_1+R(l_1),l_2+R(l_2)]\in L_x$ for any
$l_1,l_2\in\goth{h}_x^\perp$.

Choose any $l_1,l_2,l_3\in\goth{h}_x^\perp$. We are going to check
that
\begin{gather*}
Q([l_1+R(l_1),l_2+R(l_2)],l_3+R(l_3))=\\ \bigl<l_1\otimes
l_2\otimes
l_3,-\varphi+\CYB(r)-\frac{1}{2}\Alt(\delta\otimes\id)r\bigr>.
\end{gather*}
Indeed,
\begin{gather*}
\bigl<l_1\otimes l_2\otimes l_3,[r^{12},r^{13}]\bigr>=
\bigl<l_1\otimes l_2\otimes l_3,\sum_{i,j}[r_i',r_j']\otimes
r_i''\otimes r_j''\bigr>=\\
\bigl<l_1,\sum_{i,j}[\bigl<l_2,r_i''\bigr>r_i',\bigl<l_3,r_j''\bigr>r_j']\bigr>=
Q(l_1,[R(l_2),R(l_3)])=\\
Q([l_1,R(l_2)],R(l_3)).
\end{gather*}
Similarly,
\begin{equation*}
\bigl<l_1\otimes l_2\otimes l_3,[r^{12},r^{23}]\bigr>=
Q([R(l_1),l_2],R(l_3)),
\end{equation*}
\begin{equation*}
\bigl<l_1\otimes l_2\otimes l_3,[r^{13},r^{23}]\bigr>=
Q([R(l_1),R(l_2)],l_3).
\end{equation*}
It is easy to see that $\frac{1}{2}\Alt(\delta\otimes\id)r=
(\delta\otimes\id)r+ \tau(\delta\otimes\id)r+
\tau^2(\delta\otimes\id)r$, where $\tau(x\otimes y\otimes z)=
z\otimes x\otimes y$. We have
\begin{gather*}
\bigl<l_1\otimes l_2\otimes l_3,(\delta\otimes\id)r\bigr>=
\sum_i\bigl<l_1\otimes
l_2,\delta(r_i')\bigr>\bigl<l_3,r_i''\bigr>=\\
\sum_i\bigl<[l_1,l_2]_{\delta},\bigl<l_3,r_i''\bigr>r_i'\bigr>=
-Q([l_1,l_2],R(l_3)),
\end{gather*}
\begin{equation*}
\bigl<l_1\otimes l_2\otimes l_3,\tau(\delta\otimes\id)r\bigr>=
-Q([R(l_1),l_2],l_3),
\end{equation*}
\begin{equation*}
\bigl<l_1\otimes l_2\otimes l_3,\tau^2(\delta\otimes\id)r\bigr>=
-Q([l_1,R(l_2)],l_3),
\end{equation*}
\begin{equation*}
\bigl<l_1\otimes l_2\otimes l_3, \varphi\bigr>=-Q([l_1,l_2],l_3).
\end{equation*}

Adding up all the terms on the right hand side and using the fact
that $Q([R(l_1),R(l_2)],R(l_3))=0$ we see that
\begin{gather*}
Q([l_1+R(l_1),l_2+R(l_2)],l_3+R(l_3))=\\ \bigl<l_1\otimes
l_2\otimes
l_3,-\varphi+\CYB(r)-\frac{1}{2}\Alt(\delta\otimes\id)r\bigr>.
\end{gather*}

The r.h.s.\ of this equality vanishes for any
$l_1,l_2,l_3\in(\goth{g/h}_x)^*$ iff the image of
$\varphi-\CYB(r)+\frac{1}{2}\Alt(\delta\otimes\id)r$ in
$\bigwedge^3(\goth{g}/\goth{h}_x)$ vanishes.

The l.h.s.\ vanishes for any $l_1,l_2,l_3\in(\goth{g/h}_x)^*$ iff
$Q([l_1+R(l_1),l_2+R(l_2)],L_x)$ vanishes, i.e., since $L_x$ is
maximal isotropic, iff $[l_1+R(l_1),l_2+R(l_2)]\in L_x$.

This finishes the proof of the lemma.
\end{proof}

Suppose $v\in\bigwedge^2(\goth{g}/\goth{h}_x)$. Consider the
mapping $v\mapsto L_v$, where
\begin{equation*}
L_v=\{a+l\mid a\in\goth{g},\,l\in\goth{g}/\goth{h}_x,\,(l\otimes
\id)v=\overline{a}\}.
\end{equation*}
This is a bijection between
$\bigwedge^2(\goth{g}/\goth{h}_x)$ and the set of all Lagrangian
subspaces $L\subset \mathcal{D}(\goth{g})$ such that
$L\cap\goth{g}=\goth{h}_x$.

Further, there is a bijection between bivector fields $P_X$ on $X$
and smooth maps $x\mapsto L_x$ from $X$ to the set of all
Lagrangian subspaces in $\mathcal{D}(\goth{g})$ such that
$L_x\cap\goth{g}=\goth{h}_x$ for all $x\in X$.

From Lemmas \ref{equivariancelem}, \ref{phiplusetc_vanishlem} and
\ref{subalglemma} it follows that $(X,P_X)$ is a quasi-Poisson
homogeneous $G$-space iff the corresponding map $x\mapsto L_x$ is
$G$-equivariant, subalgebra-valued, and
$L_x\cap\goth{g}=\goth{h}_x$ for all $x\in X$.

This finishes the proof of Theorem \ref{mainth}.
\section{Twisting}
Let $G$ be a Lie group. Suppose $(P_G,\varphi)$ and
$(P_G',\varphi')$ are quasi-Poisson structures on $G$.
\begin{definition}[see \cite{YKS_qbialg_qpoiss}]
$(G,P_G',\varphi')$ is \emph{obtained by twisting} (by
$r\in\bigwedge^2\goth{g}$) from $(G,P_G,\varphi)$ if
\begin{equation*}
P_G'=P_G+r^\lambda-r^\rho,
\end{equation*}
\begin{equation*}
\varphi'=\varphi+\frac{1}{2}\Alt(\delta\otimes\id)r-\CYB(r).
\end{equation*}
\end{definition}

There is a similar relation on Lie quasi-bialgebras. Let
$\goth{g}$ be a Lie algebra, $(\delta,\varphi)$ and
$(\delta',\varphi')$ are Lie quasi-bialgebra structures on
$\goth{g}$.
\begin{definition}[see \cite{D_quasihopf_alg, YKS_qbialg_qpoiss}]
$(\goth{g},\delta',\varphi')$ is \emph{obtained by twisting} (by
$r\in\bigwedge^2\goth{g}$) from $(\goth{g},\delta,\varphi)$ if
\begin{equation*}
\delta'(x)=\delta(x)+\ad_xr\quad\text{for all $x\in\goth{g}$},
\end{equation*}
\begin{equation*}
\varphi'=\varphi+\frac{1}{2}\Alt(\delta\otimes\id)r-\CYB(r).
\end{equation*}
\end{definition}
Twisting is an equivalence relation.

If $(G,P_G',\varphi')$ is obtained by twisting from
$(G,P_G,\varphi)$ then the corresponding Lie quasi-bialgebra
$(\goth{g},\delta',\varphi')$ is obtained by twisting from
$(\goth{g},\delta,\varphi)$. The converse holds if $G$ is
connected.

Denote by $\mathcal{D}(\goth{g},\delta,\varphi)$ and
$\mathcal{D}(\goth{g},\delta',\varphi')$ the double Lie algebras
of Lie quasi-bialgebras $(\goth{g},\delta,\varphi)$ and
$(\goth{g},\delta',\varphi')$ respectively. The following result
is obtained in \cite{D_quasihopf_alg}.
\begin{theorem}\label{twist_equivth}
$(\goth{g},\delta',\varphi')$ is obtained by twisting from
$(\goth{g},\delta,\varphi)$ if and only if there exists a Lie
algebra isomorphism
$f_r:\mathcal{D}(\goth{g},\delta,\varphi)\to\mathcal{D}(\goth{g},\delta',\varphi')$
fixing all the elements of $\goth{g}$ and preserving the canonical
bilinear forms on the doubles.
\end{theorem}

Suppose that $(G,P_G',\varphi')$ is obtained by twisting from
$(G,P_G,\varphi)$. Let $r\in\bigwedge^2\goth{g}$ be the
corresponding bivector. Then
$f_r:\mathcal{D}(\goth{g},\delta,\varphi)\to\mathcal{D}(\goth{g},\delta',\varphi')$,
$f_r(a+l)=a+l+(l\otimes\id)r$ is the corresponding Lie algebra
isomorphism.

Using $f_r$ we can identify $\mathcal{D}(\goth{g},\delta,\varphi)$
and $\mathcal{D}(\goth{g},\delta',\varphi')$. Since $f_r$
preserves the canonical bilinear forms, the sets of Lagrangian
subalgebras under this identification are the same.
\begin{theorem}
Let $(X,P_X)$ be a homogeneous quasi-Poisson
$(G,P_G,\varphi)$-space. Then $(X,P_X-r_X)$ is a homogeneous
quasi-Poisson $(G,P_G',\varphi')$-space, and the map $P_X\mapsto
P_X-r_X$ is a bijection between the set of all $(G,P_G,\varphi)$-
and $(G,P_G',\varphi')$-quasi-Poisson structures on $X$.
\end{theorem}
\begin{proof}
Denote by $\Lambda$ (resp. $\Lambda'$) the set of all Lagrangian
Lie subalgebras in $\mathcal{D}(\goth{g},\delta,\varphi)$ (resp.
$\mathcal{D}(\goth{g},\delta',\varphi')$).

Theorem \ref{mainth} gives us the $G$-equivariant map $x\mapsto
L_x$ from $X$ to $\Lambda$ such that $L_x\cap\goth{g}=\goth{h}_x$
defined by
\begin{equation*}
L_x=\{a+l\mid a\in\goth{g},\ l\in\goth{h}_x^\perp,\
(l\otimes\id)P_X(x)=\overline{a}\}.
\end{equation*}
On the other hand, consider the map $x\mapsto L_x'$ from $X$ to
the set of subspaces in $\mathcal{D}(\goth{g},\delta',\varphi')$
corresponding to $P_X-r_X$:
\begin{equation*}
L_x'=\{a+l\mid a\in\goth{g},\ l\in\goth{h}_x^\perp,\
(l\otimes\id)(P_X(x)-r_X)=\overline{a}\}.
\end{equation*}

It is easy to see that $f_r(L_x)=L_x'$. Since $f_r$ is a Lie
algebra isomorphism, preserves the canonical bilinear forms on the
doubles and commutes with the action of $G$ on the doubles, the
map $x\mapsto L_x'$ is a $G$-equivariant map from $X$ to
$\Lambda'$. Since $f_r$ fixes all the points of $\goth{g}$, we
have $L_x'\cap\goth{g}=\goth{h}_x$. From Theorem \ref{mainth} it
follows that $P_X-r_X$ defines a $(G,P_G',\varphi')$-quasi-Poisson
structure on $X$.

Obviously, the map $P_X\mapsto P_X-r_X$ from the set of all
$(G,P_G,\varphi)$-quasi-Poisson structures on $X$ to the set of
all $(G,P_G',\varphi')$-quasi-Poisson structures on $X$ is
injective. Similarly, the map $P_X'\mapsto P_X'+r_X$ transforms a
$(G,P_G',\varphi')$-structure to a $(G,P_G,\varphi)$-structure.
Thus, we have a bijection.
\end{proof}
\section{Examples}
Recall that if $(G,P_G)$ is a Poisson Lie group, then the
homogeneous $G$-spaces $X=\{x\}$ and $Y=G$ admit the structure of
Poisson homogeneous $(G,P_G)$-spaces. Here we consider the
quasi-Poisson case.
\begin{example}
Let $(G,P_G,\varphi)$ be a quasi-Poisson Lie group, $X=\{x\}$ is a
ho\-mo\-ge\-neous $G$-space, $P_X=0$ is the only bivector field on
$X$. Then the (trivial) action of $G$ on $X$ is quasi-Poisson. The
corresponding Lagrangian subalgebra is $\goth{g}$.
\end{example}
\begin{example}
Consider the action of a connected quasi-Poisson Lie group
$(G,P_G,\varphi)$ on $Y=G$ by left translations. By Theorem
\ref{mainth}, there is a bijection between the set of
$G$-quasi-Poisson structures on $Y$ and the set of $G$-conjugacy
classes of Lagrangian subalgebras $L\subset\mathcal{D}(\goth{g})$
such that $L\cap\goth{g}=0$.

The map $r\mapsto
L_r=\{a+l\in\mathcal{D}(\goth{g})\mid(l\otimes\id)r=a\}$ from
$\bigwedge^2\goth{g}$ to the set of Lagrangian subspaces in
$\mathcal{D}(\goth{g})$ transversal to $\goth{g}$ is a bijection.
On the other hand, $L_r$ is a Lie subalgebra iff
$\varphi+\frac{1}{2}\Alt(\delta\otimes\id)r-\CYB(r)=0$.

Thus $Y$ can be a quasi-Poisson homogeneous $G$-space if and only
if $G$ is obtained by twisting from a Poisson Lie group. In this
case there is a 1-1 correspondence between the solutions of the
equation
\begin{equation*}
\CYB(r)-\frac{1}{2}\Alt(\delta\otimes\id)r=\varphi
\end{equation*}
and $(G,P_G,\varphi)$-quasi-Poisson structures on $Y$ given by
$P_Y=P_G+r^\lambda$.
\end{example}
Let us also introduce the following purely quasi-Poisson example.
\begin{example}
Suppose $\goth{g}$ is a finite-dimensional Lie algebra with a
non-degenerate invariant symmetric bilinear form $(\;|\;)$. Let
$G$ be a connected Lie group such that $\Lie G=\goth{g}$. Consider
the ``Manin quasi-triple'' (see \cite{AAYKS_mpairs_mmaps})
$(\goth{a},\goth{a}_1,\goth{a}_2)$, where
$\goth{a}=\goth{g}\times\goth{g}$,
\begin{equation*}\goth{a}_1=\{(x,x)\mid
x\in\goth{g}\}\simeq\goth{g},\ \goth{a}_2=\{(x,-x)\mid
x\in\goth{g}\},
\end{equation*}
and $\goth{a}$ is equipped with a
non-degenerate invariant symmetric bilinear form
$((a,b),(c,d))\mapsto\frac{1}{2}\left((a|c)-(b|d)\right)$. It is
easy to calculate that the corresponding Lie quasi-bialgebra
structure on $\goth{g}$ is given by $\delta=0$,
$\varphi=[\Omega^{12},\Omega^{23}]=-\CYB(\Omega)$, where
$\Omega\in(S^2\goth{g})^\goth{g}$ corresponds to $(\;|\;)$. This
Lie quasi-bialgebra gives rise to the quasi-Poisson Lie group
$(G,0,\varphi)$.

Pick any $g\in G$, and consider the Lagrangian subalgebra
\begin{equation*}
L_g=\{(x,y)\mid y=\Ad_gx\}\subset\goth{a}.
\end{equation*}
It can be
shown that it corresponds to the quasi-Poisson homogeneous space
$(C_g, P)$, where $C_g\subset G$ is the conjugacy class of $g$,
and
\begin{equation*}
P(g)=(r_g\otimes l_g - l_g\otimes r_g)(\Omega).
\end{equation*}
Moreover, one
can show that $(G, P)$ is a quasi-Poisson $G$-manifold with
respect to the action by conjugation, and $(C_g, P)$ are
``quasi-Poisson $G$-submanifolds'' of $(G, P)$ (see
\cite{AAYKSEM}, where this example was introduced and studied for
a compact Lie group $G$).
\end{example}

%
\noindent E.K.: Department of Mathematics, Kharkov National
University,\\ 4 Svobody Sq., Kharkov, 61077, Ukraine; \\ Institute
for Low Temperature Physics \& Engineering,\\ 47 Lenin Avenue,
Kharkov, 61103, Ukraine\\ e-mail: {\small \tt
eugene.a.karolinsky@univer.kharkov.ua; karol@sky.net.ua}

\medskip

\noindent K.M.: Department of Mathematics, Kharkov National
University,\\ 4 Svobody Sq., Kharkov, 61077, Ukraine\\ e-mail:
{\small \tt ono@ukr.net}

\end{document}